\documentclass{IEEEtran4PSCC}
\usepackage{algorithm, algpseudocode}
\usepackage{amsfonts, amsthm}
\usepackage{xcolor}
\usepackage{caption, subcaption, tabularx, bm}

\newcommand{\lambdab}{\bm{\lambda}}
\newcommand{\mub}{\bm{\mu}} 
\newcommand{\xb}{\mathbf{x}}

\newcommand{\cb}{\mathbf{c}}
\newcommand{\db}{\mathbf{d}}
\newcommand{\bb}{\mathbf{b}}
\newcommand{\Ab}{\mathbf{A}}
\newcommand{\Bb}{\mathbf{B}}
\newcommand{\Qb}{\mathbf{Q}}
\newcommand{\Ib}{\mathbf{I}}
\newcommand{\xbl}{\underline{\mathbf{x}}}
\newcommand{\xbu}{\overline{\mathbf{x}}}
\newcommand{\Nc}{\mathcal{N}}
\newcommand{\Gc}{\mathcal{G}}
\newcommand{\Pc}{\mathcal{P}}
\newcommand{\Dc}{\mathcal{D}}
\newcommand{\Yc}{\mathcal{Y}}
\newcommand{\Ec}{\mathcal{E}}
\newcommand{\Lc}{\mathcal{L}}
\newcommand{\Ic}{\mathcal{I}}
\newcommand{\pg}{p^{\text{g}}}
\newcommand{\qg}{q^{\text{g}}}
\newcommand{\pb}{p^{\text{b}}}
\newcommand{\qb}{q^{\text{b}}}
\newcommand{\pd}{p^{\text{d}}}
\newcommand{\qd}{q^{\text{d}}}
\newcommand{\pgl}{\underline{p}^{\text{g}}}
\newcommand{\pgu}{\overline{p}^{\text{g}}}
\newcommand{\qgl}{\underline{q}^{\text{g}}}
\newcommand{\qgu}{\overline{q}^{\text{g}}}
\newcommand{\pl}{\underline{p}}
\newcommand{\pu}{\overline{p}}
\newcommand{\ql}{\underline{q}}
\newcommand{\qu}{\overline{q}}
\newcommand{\wl}{\underline{w}}
\newcommand{\wu}{\overline{w}}
\newcommand{\thetal}{\underline{\theta}}
\newcommand{\thetau}{\overline{\theta}}
\newcommand{\gsh}{g^{\text{sh}}}
\newcommand{\bsh}{b^{\text{sh}}}
\newcommand{\gs}{g^{\text{s}}}
\newcommand{\bs}{b^{\text{s}}}
\newcommand{\Mp}{M^{\text{p}}}
\newcommand{\Mq}{M^{\text{q}}}

\ifCLASSINFOpdf
   \usepackage[pdftex]{graphicx}
\else
   \usepackage[dvips]{graphicx}
\fi
\usepackage[cmex10]{amsmath}
\makeatletter
\let\old@ps@headings\ps@headings
\let\old@ps@IEEEtitlepagestyle\ps@IEEEtitlepagestyle
\def\psccfooter#1{%
    \def\ps@headings{%
        \old@ps@headings%
        \def\@oddfoot{\strut\hfill#1\hfill\strut}%
        \def\@evenfoot{\strut\hfill#1\hfill\strut}%
    }%
    \def\ps@IEEEtitlepagestyle{%
        \old@ps@IEEEtitlepagestyle%
        \def\@oddfoot{\strut\hfill#1\hfill\strut}%
        \def\@evenfoot{\strut\hfill#1\hfill\strut}%
    }%
    \ps@headings%
}
\makeatother

\psccfooter{%
        \parbox{\textwidth}{\hrulefill \\ \small{23rd Power Systems Computation Conference} \hfill \begin{minipage}{0.2\textwidth}\centering \vspace*{4pt} \includegraphics[scale=0.06]{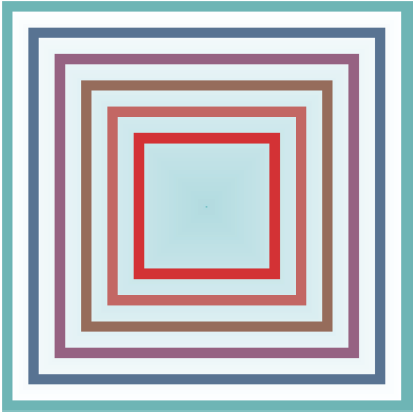}\\\small{PSCC 2024} \end{minipage} \hfill \small{Paris-Saclay, France --- June 4 -- 7, 2024}}%
}

\begin{document}
%
\title{A GPU-based Distributed Algorithm for Linearized Optimal Power Flow in Distribution Systems}


\author{
\IEEEauthorblockN{Minseok Ryu and Geunyeong Byeon}
\IEEEauthorblockA{School of Computing and Augmented Intelligence \\
Arizona State University\\
Tempe, AZ, USA 
}
\and
\IEEEauthorblockN{Kibaek Kim}
\IEEEauthorblockA{Mathematics and Computer Science Division \\
Argonne National Laboratory\\
Lemont, IL, USA 
}
}


\maketitle

\begin{abstract}
We propose a GPU-based distributed optimization algorithm, aimed at controlling optimal power flow in multi-phase and unbalanced distribution systems. Typically, conventional distributed optimization algorithms employed in such scenarios rely on parallel execution with multiple CPUs. However, this often leads to significant computation time primarily due to the need for optimization solvers to solve subproblems for every iteration of the algorithms. To address this computational challenge, we propose a distributed optimization algorithm that eliminates solver dependencies and harnesses GPU acceleration. 
The central idea involves decomposing networks to yield subproblems with closed-form solutions based on matrix operations that GPU can efficiently handle. We demonstrate the computational performance of our approach through numerical experiments on four IEEE test instances ranging from 13 to 8500 buses. Our results affirm the scalability and superior speed of our GPU-based approach compared to the CPU-based counterpart. 
\end{abstract}


\thanksto{\noindent Submitted to the 23rd Power Systems Computation Conference (PSCC 2024).}
\thanksto{\noindent This material is based upon work supported by the U.S. Department of Energy, Office of Science, under contract number DE-AC02-06CH11357. We gratefully acknowledge the computing resources provided on Bebop and Swing, a high-performance computing cluster operated by the Laboratory Computing Resource Center at Argonne National Laboratory. This research was conducted when the first author was affiliated with Argonne National Laboratory.}

\section{Introduction} \label{sec:intro} 
The modern evolution of distribution systems, exemplified by the rising integration of renewable and distributed energy sources such as solar photovoltaic, energy storage systems, and electric vehicles, is imposing a greater need for advanced operational and planning strategies based on distributed Optimal Power Flow (OPF).
The significance of distributed OPF solution within distribution systems is expanding due to its inherent benefits in terms of scalability, adaptability, privacy, and robustness in comparison to centralized optimization \cite{patari2021distributed}.

Distributed OPF in both distribution and transmission systems has been widely studied in the literature (e.g., \cite{Kim_2000, sun2013fully, peng2014distributed, Erseghe_2014,  mhanna2018adaptive, muhlpfordt2021distributed, ryu2021privacy, alkhraijah2022assessing, ryu2023differentially, Inaolaji_2023}).
It is often formulated as the following consensus optimization problem:
\begin{subequations}
\label{DO_model}
\begin{align}
\min \ &  \sum_{s=1}^S f_s(\xb_s) \label{DO_model_obj} \\
\mbox{s.t.} \ & \xb_s \in \mathcal{X}_s, \ \forall s \in [S], \label{DO_model_const} \\
& \sum_{s=1}^S \Bb_s \xb_s = \mathbf{0}, \label{DO_model_consensus}
\end{align}
\end{subequations}
where $\xb_s$, $\mathcal{X}_s$, and $f_s$ represent variables, feasible region, and objective function, respectively, defined for each subsystem $s \in [S]$ of the entire power network.
Eq. \eqref{DO_model_consensus} refers to consensus constraints that establish interconnections among subsystems within the network, thereby guaranteeing the feasibility of the solutions acquired. 
Various distributed optimization algorithms have been proposed for solving \eqref{DO_model}, contingent upon the chosen network decomposition strategy, resulting in different $S$, and the specific modeling of OPF, resulting in different $g_s$.
The widely used algorithm is the alternating direction method of multipliers (ADMM) which is an iterative algorithm that solves the augmented Lagrangian formulation of \eqref{DO_model}, derived by penalizing \eqref{DO_model_consensus} in the objective function. 
At each iteration of ADMM, it is possible to simultaneously solve subproblems, each necessitating the utilization of optimization solvers running on CPUs. 

In a recent study \cite{kim2022accelerated}, the authors have made significant improvement on computing distributed OPF in transmission systems by employing the computational power of multiple GPUs. 
Specifically, they decompose the power network componentwisely as in \cite{mhanna2018adaptive}, resulting in a bunch of small subproblems, and apply the two level ADMM algorithm \cite{sun2021two}.
While a significant portion of the subproblems admits closed-form solutions, a subset of these subproblems takes the form of bound-constrained nonlinear programming (NLP) models. 
To solve the models, a trust-region Newton algorithm \cite{kim2021leveraging} is employed, which fully utilizes multiple GPUs.

Motivated by this achievement, this paper aims to utilize multiple GPUs to accelerate the computation of distributed OPF in distribution systems which are more complicated compared with transmission system mainly due to the existence of multi-phase and unbalanced power loads connected by either wye or delta configuration.
Specifically, we employ the linearization technique developed for computing OPF solution in the unbalanced distribution systems \cite{byeon2022linear} and propose a distributed algorithm that is free from optimization solvers because each subproblem to solve for each iteration of the algorithm admits a closed-form solution expression based on matrix operations, leading to extremely fast computation.
Notably, the closed-form solution expression holds for any size of $S \geq 1$ in \eqref{DO_model}.
The proposed algorithm is different from existing GPU-based optimization solvers (e.g., cuOSQP \cite{schubiger2020gpu}) that is for $S=1$ in \eqref{DO_model} while our algorithm can be used for $S \geq 1$. 
Our contributions are summarized as follows:
\begin{enumerate}
  \item development of a distributed optimization algorithm that does not rely on optimization solvers due to the existence of closed form solution expressions; 
  \item implementation of the algorithm that can efficiently run on NVIDIA GPU architecture, as well as CPUs; and  
  \item demonstration of the greater computational performance of the proposed GPU-based algorithm compared to the multi-CPU counterpart.
\end{enumerate}

The remainder of this paper is organized as follows. 
In Section \ref{sec:math_form}, we describe the mathematical formulation for the linearized distributed OPF in distribution systems with unbalanced power loads and present existing decomposition methods for solving the model.
Then we present the proposed GPU-based distributed optimization algorithm in Section \ref{sec:proposed}.
Finally, in Section \ref{sec:experiments}, we numerically showcase the superiority of our approach with respect to computation time.

\begin{table} 
\centering
\caption{Nomenclature}  
\begin{tabularx}{0.5\textwidth}{ll}      
\hline
\multicolumn{2}{c}{Sets} \\ \hline
  $\Nc$, $\Gc$, $\Ec$ & a set of buses, generators, and lines, respectively\\  
  $\Lc_i$, $\Yc_i$, $\Dc_i$ & a set of loads, wye loads, delta loads at bus $i \in \Nc$ \\
  $\Pc_i \subseteq \{1,2,3\}$ & a set of phases at component $i \in \Nc \cup \Gc \cup \Ec \cup \Lc$ \\
  \hline 
  \multicolumn{2}{c}{Parameters} \\ \hline
  $\pgl_{k\phi} , \pgu_{k\phi}, \qgl_{k\phi} , \qgu_{k\phi}$ & bounds on the power generation at $k \in \Gc$ and $\phi \in \Pc_k$.  \\
  $\wl_{i\phi}, \wu_{i\phi}$ & bounds on the squared magnitude of voltage  \\ & at $i \in \Nc$ and $\phi \in \Pc_i$. \\
  $\gsh_{i \phi}, \bsh_{i \phi}$ & shunt conductance and shunt susceptance \\ & 
  at $i \in \Nc$ and $\phi \in \Pc_i$. \\  
  $r_{e \phi \phi'}, x_{e \phi \phi'}$ & resistance and reactance at $e \in \Ec$ and $\phi, \phi' \in \Pc_{e}$. \\
  $\thetal_{e}, \thetau_{e}$ & bounds on the angle difference at line $e \in \Lc$.\\
  $g^{\text{s}}_{e i j \phi}, b^{\text{s}}_{e ij \phi }$ & shunt conductance and shunt susceptance \\ & at from-bus $i$ of $e \in \Ec$ on phase $\phi \in \Pc_{e}$. \\
  $\tau_{e \phi}$ & tap ratio at $e \in \Ec$ and $\phi \in \Pc_{e}$. \\       
\hline
\multicolumn{2}{c}{Variables} \\ \hline
$\pg_{k\phi}, \qg_{k\phi}$ & power generation at $k \in \Gc$ and $\phi \in \Pc_k$.\\
$w_{i\phi} $ & squared magnitude of voltage at $i \in \Nc$ and $\phi \in \Pc_i$. \\
$\pb_{l\phi},  \qb_{l\phi}, \pd_{l\phi}, \qd_{l\phi}$ & voltage dependent load at $l \in \Lc$ and $\phi \in \Pc_l$.  \\
$p_{e ij \phi}, q_{e ij \phi}$ & power flow on line $e_{ij} \in \Ec$ and phase $\phi \in \Pc_{e}$. \\ \hline
\end{tabularx}  
\label{Table:Nomenclature}
\end{table}
 
\section{Mathematical Formulations} \label{sec:math_form}
In this section we formulate a linearized OPF model for a distribution system decomposed into several subsystems.  

\subsection{Linearized OPF model for distribution systems}
Based on the notations in Table \ref{Table:Nomenclature}, we present a centralized LP model proposed in \cite{byeon2022linear} without derivation. 
\subsubsection{Operational bound constraints} We consider
\begin{subequations}
\label{operational}
\begin{align}
& \pg_{k\phi} \in [\pgl_{k\phi},\pgu_{k\phi}], \qg_{k\phi} \in [\qgl_{k\phi},\qgu_{k\phi}], \ \forall k \in \Gc, \phi \in \Pc_k, \\
& w_{i\phi} \in [\wl_{i\phi}, \wu_{i\phi}], \ \forall i \in \Nc,  \phi \in \Pc_i, \\
& p_{eij\phi},p_{eji\phi} \in [\pl_{eij\phi},\pu_{eij\phi}], \ \forall (e,i,j) \in \Ec, \phi \in \Pc_e, \\
& q_{eij\phi},q_{eji\phi} \in [\ql_{eij\phi},\qu_{eij\phi}], \ \forall (e,i,j) \in \Ec, \phi \in \Pc_e ,
\end{align}
\end{subequations}
to ensure that power generation, voltage magnitude, and power flow, respectively, are within certain bounds.

\subsubsection{Power balance equations}
For every bus $i \in \Nc$ and phase $\phi \in \Pc_i$, the power balance equations are given by
\begin{subequations}
\label{balance_eqns}
\begin{align}
& \sum_{ (e, i,j) \in \Ec_i } p_{e i j  \phi} + \sum_{l \in \Lc_i} \pb_{l \phi} + \gsh_{i \phi} w_{i \phi} = \sum_{ k \in \Gc_i } \pg_{k \phi}, \\
& \sum_{ (e, i,j) \in \Ec_i } q_{e i j \phi} + \sum_{l \in \Lc_i} \qb_{l \phi}  - \bsh_{i \phi} w_{i \phi} = \sum_{ k \in \Gc_i } \qg_{k \phi},  
\end{align}  
\end{subequations}
which ensure power balance at every bus in the network.
Note that $\pb_{l\phi}$ and $\qb_{l\phi}$ represent real and reactive power loads on phase $\phi$ drawn by a load $l$ from the bus to which it is attached. 

\subsubsection{Voltage dependent load model}
For every bus $i \in \Nc$, the voltage dependent load model in \cite{byeon2022linear} is given by  
\begin{subequations}
\label{volt_load_model}  
\begin{align}
& \pd_{l\phi} = \frac{a_{l\phi} \alpha_{l\phi}}{2} (\widehat{w}_{l\phi}-1)+ a_{l\phi}, \ \forall l \in \Lc_i, \forall \phi \in \Pc_i, \label{VDLM-1}\\
& \qd_{l\phi} = \frac{b_{l\phi} \beta_{l\phi}}{2} (\widehat{w}_{l\phi}-1) + b_{l\phi},  \   \forall l \in \Lc_i, \forall \phi \in \Pc_i, \label{VDLM-2}\\
&  \widehat{w}_{l \phi} = w_{i\phi}, \  \forall l \in \Yc_i, \forall \phi \in   \Pc_i, \label{VDLM-3}\\
&  \widehat{w}_{l \phi} = 3 w_{i\phi}, \ \forall l \in \Dc_i, \forall \phi \in \Pc_i, \label{VDLM-4}
\end{align}
where, for fixed $l \in \Lc_i$ and $\phi \in \Pc_i$,  
$\pd_{l\phi}$ and $\qd_{l\phi}$ are real and reactive power consumption,  
$\alpha_{l\phi}$ and $\beta_{l\phi}$ are nonnegative scalars given as inputs that determines the load type, 
$a_{l\phi}$ and $b_{l\phi}$ are also given input data determined by reference values for real power load, reactive power load, and voltage magnitude applied to the load, 
$\widehat{w}_{l\phi}$ in \eqref{VDLM-1} and \eqref{VDLM-2} are local variables determined by \eqref{VDLM-3} or \eqref{VDLM-4} relying on the configurtation of $l$.

For the wye connected load $l \in \Yc$ on phase $\phi \in \Pc_l$, the relationship between $\pb$ and $\pd$ is exact:
\begin{align}
\pb_{l \phi} = \pd_{l \phi}, \ \ \qb_{l \phi} = \qd_{l \phi}. \label{VDLM-5}
\end{align}
For the delta connected load, $l \in \Dc$, the relationship between $\pb$ and $\pd$ is given by 
\begin{align} 
& \sum_{\phi \in \Pc_l} \pb_{l\phi} - \pd_{l\phi} = \sum_{\phi \in \Pc_l} \qb_{l\phi} - \qd_{l\phi} = 0,  \label{VDLM-6}\\ 
&  \frac{3}{2} \pb_{l2} - \frac{\sqrt{3}}{2} \qb_{l2}  = {\pd_{l2} + \frac{1}{2} \pd_{l1} - \frac{\sqrt{3}}{2} \qd_{l1}},  \label{VDLM-7}\\
& \frac{\sqrt{3}}{2} \pb_{l2} + \frac{3}{2} \qb_{l2} = {\frac{\sqrt{3}}{2} \pd_{l1} +\frac{1}{2} \qd_{l1} + \qd_{l2}},  \label{VDLM-8}\\
& \sqrt{3} \qb_{l2} + \frac{3}{2} \pb_{l3}-\frac{\sqrt{3}}{2} q^{b}_{l3} = {\frac{1}{2} \pd_{l1} + \frac{\sqrt{3}}{2} \qd_{l1} + \pd_{l3}},  \label{VDLM-9} \\
& - \sqrt{3}\pb_{l2} + \frac{\sqrt{3}}{2} \pb_{l3} + \frac{3}{2} q^{b}_{l3}   = {- \frac{\sqrt{3}}{2} \pd_{l1} + \frac{1}{2} \qd_{l1}  + \qd_{l3}} ,  \label{VDLM-10} 
\end{align}
\end{subequations}
which is derived in \cite{byeon2022linear}.

\subsubsection{Linearized power flow equations}
For every line $(e,i,j) \in \Ec$ and phase $\phi \in \Pc_e$, the linearized power flow equations are given by
\begin{subequations}
\label{power_flow_eqns}  
\begin{align} 
& p_{e ij \phi} + p_{e ji \phi} = \gs_{e i j \phi} w_{i \phi} + \gs_{e ji \phi} w_{j \phi}, \label{powerloss-1} \\
& q_{e ij \phi} + q_{e ji \phi} = -\bs_{e ij \phi}  w_{i \phi} -\bs_{e ji \phi} w_{j \phi}, \label{powerloss-2} \\ 
& w_{i \phi} = \tau_{e \phi} w_{j \phi} - \sum_{\psi \in \Pc_{e}} \Mp_{e \phi \psi} (p_{e ij \psi}-\gs_{e ij \psi} w_{i \psi}) \nonumber \\
& \ \ \ \ \ \ \ \ \ \ \ \  - \sum_{\psi \in \Pc_{e}} \Mq_{e \phi \psi} (q_{e ij \psi}+\bs_{e ij \psi} w_{i \psi}), \label{voltage_mag_diff}   
\end{align}
\end{subequations}
where $\Mp_{e \phi \psi}$ and $\Mq_{e \phi \psi}$ are $(\phi, \psi)$-th element of $\Mp_{e}$ and $\Mq_{e}$, respectively, which are defined as

\begin{align*} 
  & \Mp_{e} = 
  \begin{bmatrix}
    -2r_{\ell 11}  & r_{\ell 12} - \sqrt{3}x_{\ell 12} & r_{\ell 13} + \sqrt{3}x_{\ell 13} \\
    r_{\ell 21} + \sqrt{3}x_{\ell 21} & -2r_{\ell 22} & r_{\ell 23} - \sqrt{3}x_{\ell 23} \\ 
    r_{\ell 31} - \sqrt{3} x_{\ell 31} & r_{\ell 32} + \sqrt{3}x_{\ell 32} & -2r_{\ell 33}
  \end{bmatrix} \\
  & \Mq_{e} = 
  \begin{bmatrix}
    -2x_{\ell 11}  & x_{\ell 12} + \sqrt{3}r_{\ell 12} & x_{\ell 13} - \sqrt{3}r_{\ell 13} \\
    x_{\ell 21} - \sqrt{3} r_{\ell 21} & -2x_{\ell 22} & x_{\ell 23} + \sqrt{3}r_{\ell 23} \\ 
    x_{\ell 31} + \sqrt{3} r_{\ell 31} & x_{\ell 32} - \sqrt{3}r_{\ell 32} & -2x_{\ell 33}
  \end{bmatrix} 
\end{align*}

\subsubsection{Linearized OPF model}
To summarize, we have 
\begin{subequations}
\label{lindist3flow}  
\begin{align} 
\min \ &  \sum_{k \in \Gc} \sum_{\phi \in \Pc_k} \pg_{k \phi} \label{lindist3flow_obj}\\
\mbox{s.t.} \ 
& \eqref{operational}, \eqref{balance_eqns}, \eqref{volt_load_model}, \eqref{power_flow_eqns},  
\end{align}
\end{subequations}
which can be represented as the following abstract LP form:
\begin{subequations}
\label{LP_model}
\begin{align}
\min \ & \cb^{\top} \xb \label{LP_model_obj} \\
\mbox{s.t.} \ & \Ab \xb = \bb, \label{LP_model_eqn}\\
\ & \xbl \leq \xb \leq \xbu, \label{LP_model_bnd}
\end{align}
where
\begin{align}
\xb = 
\begin{bmatrix}
  & \{\pg_{k\phi}, \qg_{k\phi}\}_{k\in\Gc, \phi \in \Pc_k}   \\
  & \{ w_{i\phi}\}_{i\in\Nc, \phi \in \Pc_i} \\
  & \{\pb_{l\phi}, \qb_{l\phi}, \pd_{l\phi}, \qd_{l\phi} \}_{l \in \Lc_i,  \phi \in \Pc_l } \\
  & \{ p_{eij\phi}, q_{eij\phi}, p_{eji\phi}, q_{eji\phi} \}_{(e,i,j) \in \Ec, \phi \in \Pc_e}
\end{bmatrix}, \nonumber
\end{align}
\end{subequations}
\eqref{LP_model_obj} corresponds to \eqref{lindist3flow_obj},
\eqref{LP_model_eqn} corresponds to \eqref{balance_eqns}, \eqref{volt_load_model}, \eqref{power_flow_eqns}, 
and
\eqref{LP_model_bnd} corresponds to \eqref{operational}.

\subsection{Distributed OPF model} \label{sec:dist_opf}
The centralized LP model \eqref{LP_model} can be written as a distributed optimization model of the form \eqref{DO_model}:
\begin{subequations}
\label{LP_model_decomposed}
\begin{align}
\min \ &  \sum_{s=1}^S \cb^{\top}_s \xb_s \\
\mbox{s.t.} \ & \Ab_s \xb_s = \bb_s, \ \ \xb_s \in [\xbl_s, \xbu_s], \ \forall s \in [S], \\
& \sum_{s=1}^S \Bb_s \xb_s = \mathbf{0}. \label{LP_model_decomposed:consensus}
\end{align}
\end{subequations}
where $\Ab_s$, $\bb_s$, $\cb_s$, $\xbl_s$, and $\xbu_s$ are given input data that describe each subsystem $s$.
Since some variables belong to multiple subsystems, we need \eqref{LP_model_decomposed:consensus} to ensure feasibility of $\{\xb_s\}_{s=1}^S$.

Since the model \eqref{LP_model_decomposed} without the concensus constraint \eqref{LP_model_decomposed:consensus} is separable over $s \in [S]$, existing decomposition algorithms (e.g., dual decomposition or ADMM) can be utilized, which solves the subproblems in parallel for every iteration of the algorithms.
The computational intricacy associated with each subproblem primarily depends on how the network is decomposed.
For example, if the network is partitioned into a smaller number of relatively larger subnetworks, then each subproblem needs to be solved by an LP solver. 
Conversely, in a case where the network is decomposed component-wisely as in \cite{mhanna2018adaptive}, it results in a set of generator, bus, and branch subproblems, where generator and bus subproblems admit closed-form solutions, whereas the branch subproblems demand the utilization of an LP solver.

\section{The proposed distributed optimization algorithm running on GPUs} \label{sec:proposed}
In contrast to the existing approaches as depicted in Section \ref{sec:dist_opf}, the proposed approach is free from optimization solvers because each subproblem has a closed-form solution expression based on the matrix operations.
Therefore, our approach can offer enhanced scalability, particularly when GPUs are leveraged, given their well-documented efficiency in conducting large-scale matrix-vector multiplications through parallel computation within a shared memory environment.
It is noteworthy, however, that in a certain case (e.g., component-based decomposition \cite{mhanna2018adaptive}), the size of matrix $\Ab_s$ associated with a given subsystem $s \in [S]$ may be so small that CPU can also conduct the matrix-vector multiplication efficiently.
We will also discuss how to utilize GPUs in such a situation.

\subsection{Our Model}
First, we rewrite \eqref{LP_model} as
\begin{subequations}
\label{OurModel}
\begin{align}
\min \ & \cb^{\top} \xb \\
\mbox{s.t.} \ & \Ab_s \xb_s = \bb_s, \ \forall s \in [S], \label{OurModel_subconst}  \\
& \Bb_s \xb = \xb_s, \ \forall s \in [S], \label{OurModel_consensus} \\
& \xbl \leq \xb \leq \xbu, \label{OurModel_bound} 
\end{align}
\end{subequations}
where, for each subsystem $s \in [S]$, $\Ab_s \in \mathbb{R}^{m_s \times n_s}$, $\bb_s \in \mathbb{R}^{m_s}$ are given input data, and $\Bb_s \in \mathbb{R}^{n_s \times n}$ is a $0-1$ matrix wherein the row sums equate to $1$, while the column sums are either $0$ or $1$.
This matrix serves as a mapping operator, facilitating the transformation of a global variables vector $\xb \in \mathbb{R}^n$ into a local variables vector $\xb_s \in \mathbb{R}^{n_s}$. 
The main difference of \eqref{OurModel} from \eqref{LP_model_decomposed} is that bound constraints \eqref{OurModel_bound} are isolated from the local equality constraints \eqref{OurModel_subconst}, which is a key idea leading to a distributed algorithm that can run on GPUs. 

It is worth noting that the role of a system operator is reduced from computing \eqref{LP_model} in a centralized manner to managing the entire network within a given operational bounds in \eqref{OurModel}.
Concurrently, multiple agents are entrusted with the control and management of their respective subnetworks, drawing upon their own distinct sets of data.
In this distributed setting, the operator and the agents communicate iteratively until a global optimum is achieved. 
To this end, for every $t$-th iteration, the operator receives local solutions $\lambdab_s^{(t)}$ and $\xb_s^{(t)}$ from every agents $s \in [S]$ and solves the following optimization problem:
\begin{align}
\min_{\xb \in [\xbl, \xbu]} \cb^{\top}\xb + \sum_{s=1}^S \Big\{ \langle \lambdab_s^{(t)}, \Bb_s \xb  \rangle + \frac{\rho}{2} \|\Bb_s \xb - \xb_s^{(t)} \|^2 \Big\}, \label{ADMM-1} 
\end{align}
which computes an optimal solution $\xb^{(t+1)}$.
Then, each agent controlling a subsystem $s \in [S]$ receives $\xb^{(t+1)}$ from the operator and solves the following optimization problem:
\begin{align}
\min_{\xb_s \in \mathbb{R}^{n_s}} \ & - \langle \lambdab_s^{(t)}, \xb_s  \rangle + \frac{\rho}{2} \|\Bb_s \xb^{(t+1)} - \xb_s \|^2  \nonumber \\     
\mbox{s.t.} \ & \Ab_s \xb_s = \bb_s, \label{ADMM-2} 
\end{align}
and update the dual variables
\begin{align}
\lambdab^{(t+1)}_s = \lambdab^{(t)}_s + \rho (\Bb_s \xb^{(t+1)} - \xb_s^{(t+1)}).   \label{ADMM-3} 
\end{align}  
Note that \eqref{ADMM-1}--\eqref{ADMM-3} collectively constitute one iteration of the conventional ADMM algorithm for solving \eqref{OurModel}, derived from the following augmented Lagrangian formulation:
\begin{align*}
\max_{ \{\lambdab_s \in \mathbb{R}^{m_s} \}_{s} }  \min \ & \cb^{\top} \xb + \sum_{s=1}^S  \langle \lambdab_s, \Bb_s \xb - \xb_s \rangle + \frac{\rho}{2} \|\Bb_s \xb - \xb_s \|^2   \\
\mbox{s.t.} \ & \eqref{OurModel_subconst}, \eqref{OurModel_bound}.
\end{align*}
where $\lambdab_s \in \mathbb{R}^{m_s}$ is a dual variables vector associated with the consensus constraint \eqref{OurModel_consensus}.

\subsection{Closed-form solution expressions}
In this section we present closed-form solution expressions for subproblems \eqref{ADMM-1} and \eqref{ADMM-2}.

First, we observe that the subproblem \eqref{ADMM-1} is separable over each element of $\xb$. 
To see this, we introduce a set $\Ic_{si} := \{ j \in [n_s] : (\Bb_s)_{j,i} = 1 \}$ for all $s \in [S]$ and $i \in [n]$.
With this set, \eqref{ADMM-1} can be decomposed into $n$-subproblems, each $i \in [n]$ of which is given by
\begin{align}
\min_{x_i \in [\underline{x}_i, \overline{x}_i]} \ \big(\cb + \sum_{s=1}^S \Bb^{\top}_s \lambdab_s^{(t)} \big)_i  x_i + \frac{\rho}{2} \sum_{s=1}^S  \sum_{j \in \Ic_{si}} \big(x_i - (\xb^{(t)}_s)_j \big)^2, \nonumber
\end{align}
which is an 1-D optimization problem with a convex quadratic objective function and a bound constraint.
The closed-form solution of the problem is given by
\begin{align}
x_i^{(t+1)} = \min \ \{ \max \{ \widehat{x}_i, \underline{x}_i \}, \ \overline{x}_i \} \label{closed_1}
\end{align}
where
\begin{align}
\widehat{x}_i = \frac{-1}{\rho \sum_{s=1}^S | \Ic_{si}| } \Big\{  \big(\cb + \sum_{s=1}^S \Bb^{\top}_s \lambdab_s^{(t)} \big)_i - \sum_{s=1}^S \sum_{j \in \Ic_{si}} (\xb^{(t)}_s)_j   \Big\}.  \nonumber
\end{align}
Second, for each subsystem $s \in [S]$, we rewrite the subproblem \eqref{ADMM-2} as the following quadratic program:
\begin{subequations}
\begin{align}
\min \ & \frac{1}{2} \xb^{\top}_s \Qb_s \xb_s + \langle \db^{(t)}_s, \xb_s \rangle \nonumber \\
\mbox{s.t.} \ & \ \Ab_s \xb_s = \bb_s, \nonumber
\end{align}
\end{subequations}
where $\Qb_s=\rho \Ib_{n_s}$ and $\db^{(t)}_s=-\rho \Bb_s \xb^{(t+1)} -\lambdab_s^{(t)}$.
Without loss of generality, we assume that $\Ab_s$ is a full row-rank matrix.
If $\Ab_s$ is not given as a full row-rank matrix, row reduction techniques can be applied to $\Ab_s \xb_s = \bb_s$ as a preprocessing step to rectify this. 
Since $\Ab_s$ is a full row-rank matrix and $\Qb_s$ is a symmetric positive definite matrix, a closed-form solution expression can be readily derived.
To see this, we introduce dual variables $\mub_s$ associated with $\Ab_s \xb_s = \bb_s$. Then the Lagrangian function is given by
\begin{align}
\max_{\mub \in \mathbb{R}^{m_s} } \min_{\xb_s \in \mathbb{R}^{n_s}} \frac{1}{2} \xb^{\top}_s \Qb_s \xb_s + \langle \db^{(t)}_s, \xb_s \rangle  + \langle \mub_s , \Ab_s \xb_s - \bb_s \rangle.\nonumber
\end{align}
For given $\mub_s \in \mathbb{R}^{m_s}$, one can derive from the first optimality condition that
\begin{align}
\xb_s = \frac{1}{\rho}   (-\db^{(t)}_s - \Ab^{\top}_s\mub_s). \nonumber
\end{align}
Plugging $\xb_s$ into $\Ab_s \xb_s = \bb_s$ yields
\begin{align}
\mub = - \rho (\Ab_s \Ab^{\top}_s)^{-1} \bb_s -  (\Ab_s \Ab^{\top}_s)^{-1} \Ab_s \db^{(t)}_s.\nonumber
\end{align}
Therefore, we have
\begin{subequations}
\label{closed_2}
\begin{align}
\xb^{(t+1)}_s = \frac{1}{\rho} \overline{\Ab}_s \db^{(t)}_s + \overline{\bb}_s
\end{align}
where
\begin{align}
& \overline{\Ab}_s := \Ab^{\top}_s (\Ab_s \Ab^{\top}_s)^{-1} \Ab_s - \Ib_{n_s} \\
& \overline{\bb}_s := \Ab^{\top}_s (\Ab_s \Ab^{\top}_s)^{-1} \bb_s
\end{align}
\end{subequations}
Note that both $\overline{\Ab}_s$ and $\overline{\bb}_s$ are computed only once at a preprocessing step.


\subsection{Termination Criterion}
We use the standard primal and dual residuals to determine the termination of Algorithm \ref{alg:admm}, which are given as
\begin{align}
 \text{pres}^{(t)} \leq  \epsilon_{\text{prim}}^{(t)}, \ \   \text{dres}^{(t)} \leq \epsilon_{\text{dual}}^{(t)}, \label{termination}
\end{align}
where
\begin{align*}
& \textstyle \text{pres}^{(t)} := \sqrt{\sum_{s=1}^S \| \mathbf{B}_s \xb^{(t)} - \xb_s^{(t)}\|^2}  , \\
& \textstyle  \text{dres}^{(t)} :=  \rho \sqrt{ \sum_{s=1}^S \| \mathbf{B}^{\top}_s(\xb^{(t)}_s - \xb^{(t-1)}_s) \|^2 } , \\
& \textstyle \epsilon_{\text{prim}}^{(t)}:= \epsilon_{\text{rel}} \max \Big\{  \sqrt{\sum_{s=1}^S \|\Bb_s \xb^{(t)} }\|^2, \ \sqrt{\sum_{s=1}^S \| \xb^{(t)}_s \|^2 } \Big \}, \\
& \textstyle \epsilon_{\text{dual}}^{(t)}:= \epsilon_{\text{rel}}     \sqrt{\sum_{s=1}^S \|\Bb_s^{\top} \lambdab^{(t)}_s }\|^2 , \\
\end{align*} 
and $\epsilon_{\text{rel}}$ is a tolerance level (e.g., $10^{-2}$ or $10^{-3}$).

\subsection{The proposed distributed algorithm}
The proposed algorithm is summarized in Algorithm \ref{alg:admm}.
After initializing dual and local solutions as in line 1, we precompute the matrices and vectors used for local updates at each subsystem $s \in [S]$ as in lines 2-3.
Then, the proposed iterative algorithm runs as in lines 4-10 until a termination criterion \eqref{termination} is satisfied.
For every iteration $t$ of the algorithm, the computational steps are taken to update a global solution as in line 5, local solutions as in line 7, and dual solutions as in line 8.

\begin{algorithm}[h!]
\caption{The proposed distributed optimization algorithm }\label{alg:admm}
\begin{algorithmic}[1]
\State Initialization: $t \gets 1$, $\lambdab^{(t)}_s$, $\xb^{(t)}_s$, and $\xb^{(t,E+1)} \gets \xb^{(t)}_s$ for all $s \in [S]$
\Statex Precomputation:
\For { $s \in [S]$ }{ in parallel}
\begin{align*}
  & \overline{\Ab}_s := \Ab^{\top}_s (\Ab_s \Ab^{\top}_s )^{-1}  \Ab_s - \Ib_{n_s} \\
  & \overline{\bb}_s := \Ab^{\top}_s (\Ab_s \Ab^{\top}_s )^{-1}\bb_s
\end{align*}
\EndFor
\While{\eqref{termination} is not satisfied}
\State Compute $\xb^{(t+1)}$ by \eqref{closed_1}
\For { $s \in [S]$ }{ in parallel}
\State Compute $\xb^{(t+1)}_s$ by \eqref{closed_2}
\State Compute $\lambdab^{(t+1)}_s$ by \eqref{ADMM-3}
\EndFor 
\EndWhile 
\end{algorithmic} 
\end{algorithm}

The proposed algorithm is an ADMM algorithm composed of subproblems that admit closed-form solution expressions.
It is known that the ADMM algorithm converges to an optimal solution with sublinear rate \cite{he20121}.
A few acceleration schemes exist that could reduce the total number of iterations, including residual balancing \cite{wohlberg2017admm} and multiple local updates \cite{ryu2022differentially}.
In this paper, instead of proposing new acceleration schemes to reduce the total number of iterations, we focus on improving computation time per iteration by running the algorithm in GPUs. We remark that any existing acceleration schemes can be utilized.  

\section{Implementation Details for the GPU computations} \label{sec:implementation}
We implement the proposed distributed algorithm in Julia with the CUDA.jl package \cite{besard2018effective} for the NVIDIA GPU programming.
A GPU has been used for conducting large matrix-vector multiplications. This is mainly because of its ability of conducting such operations in parallel using multiple threads under a shared memory.
For example, when conducing $\Ab \xb$ where $\Ab \in \mathbb{R}^{m \times n}$ and $\xb \in \mathbb{R}^{n}$, each thread of a GPU conducts a  $\mathbf{a}^{\top}_i \xb$ where $\mathbf{a}^{\top}_i$ is $i$-th row of $\Ab$.
Since Algorithm \ref{alg:admm} is composed of a set of matrix-vector operations as in \eqref{closed_1}, \eqref{closed_2}, and \eqref{ADMM-3}, GPU can be efficient when each subsystem $s \in [S]$ has large matrix $\Ab_s$.
In this case, one can simply use \texttt{CuArray} supported by the CUDA.jl package for writing matrices and vectors and conduct matrix operations as done in the CPU programming.

However, in the case of subsystems with small $\Ab_s$ resulting from the component-based decomposition as in \cite{mhanna2018adaptive}, using \texttt{CuArray} will not improve the matrix-vector operations.
In this case, we construct a kernel function that assigns the matrix-vector operations $\Ab_s \xb_s$ to each thread of a GPU and enable the computation in parallel.   

For communication between multiple GPUs within a computing cluster, we use the MPI.jl package \cite{byrne2021mpi}. For every iteration of Algorithm \ref{alg:admm}, it requires a communication of intermediate solutions, including (i) $\xb^{(t+1)}$ in line 5 sent from a server to all subsystems, and (ii) $\xb^{t+1}_s$ in line 7 an $\lambdab_s^{t+1}$ in line 8 sent from each subsystem to the server.
We note that the commucation time between multiple GPUs is often higher than that between multiple CPUs because MPI requires to transfer data from GPU to CPU to communicate. 
In practice, however, multiple GPUs are located at different computing resources, thus MPI cannot be used for communication.
In this case, a remote procedure call which enables communication between multiple platforms can be utilized.
With tRPC \cite{damania2023pytorch}, for example, the communication time between multiple GPUs is about the same as that between multiple CPUs, as reported in \cite{damania2023pytorch}.
This can rule out a potential concern of communication overhead due to the introduction of GPUs.

We remark that single-precision is used for the GPU programming while double-precision is used for the CPU programming. Nonetheless, we report that it does not affect the convergence behaviour of Algorithm \ref{alg:admm} which will be reported in Section \ref{sec:experiments}.
  
\section{Numerical Results} \label{sec:experiments}
In this section we aim to demonstrate how the proposed GPU implementation of Algorithm \ref{alg:admm} accelerates the computation compared with its CPU counterpart.
Even though our algorithm can be utilized for the case $S = 1$ in \eqref{OurModel} (i.e., centralized version) and compared to the existing GPU-based solver (e.g., cuOSQP \cite{schubiger2020gpu}), it is the out-of-scope in this paper.

All numerical tests were performed on (i) Swing, a 6-node GPU cpomputing cluster (each node has 8 NVIDIA A100 40GB GPUs) and (ii) Bebop, a 1024-node CPU computing cluser (each node has 36 cores with Intel Xeon E5-2695v4 processors and 128 GB DDR4 of memory) at Argonne National Laboratory.

\subsection{Experimental Settings}
For our demonstration, we use IEEE test instances \cite{schneider2017analytic} with 13, 37, 123, and 8500 buses, which are designed for a centralized computation. 
In these instances, the three-phase loads at each bus are given as either wye or delta connected and each load is labeled as either constant power, current, or impedance load. 
These features are taken into a consideration by the LP model in Section \ref{sec:math_form}.
The size of matrix $\Ab$ in \eqref{LP_model} for each instance is reported in Table \ref{tab:size_A}.
Note that the number of columns in this table is the size of global variables $\xb$ in \eqref{OurModel}.

\begin{table}[!h]
\centering
\caption{The number of rows and columns of $\Ab$ in \eqref{LP_model} for the four test instances.}
\label{tab:size_A}
\begin{tabular}{cccc}
\hline  
IEEE13 & IEEE37& IEEE123& IEEE8500 \\ \hline
(456, 454) & (1206, 1210) & (1834, 1834) & (86114, 87285)  \\
\hline
\end{tabular}
\end{table} 

For the distributed setting, we construct $S$ subsystems in \eqref{OurModel} based on the network decomposition by its network components such as buses, transformers, and branches, as in \cite{mhanna2018adaptive}.
Specifically, we construct a graph with a set of lines, each of which represents either branch or transformer line, and a set of nodes, each of which represents either bus or a node connecting to a transformer line. 
Then we construct a set of subproblems related to each component of the graph.
Based on our observation that it is often the case that the subproblems related to leaf nodes and their connected lines are much smaller than the other subproblems, we construct a set of subsystems, each of which is constructed by combining a leaf node to its connecting line and another set of subsystems, each of which is constructed by the remaining individual node or line.  
We report the total number $S$ of subsystems we constructed in Table \ref{tab:subsystems}.

\begin{table}[!h]
  \centering
  \caption{Total number $S$ of subsystems}
  \label{tab:subsystems}
  \begin{tabular}{c|cccc}
  \hline
         & IEEE13 &  IEEE37 & IEEE123 & IEEE8500 \\ \hline    
  \# of nodes   & 29&   56& 147 & 11932 \\         
  \# of lines   & 28&   57& 146 & 14291 \\         
  \# of leaf nodes & 7 &   16& 43 & 1222 \\ 
  $S$           & 50&  97& 250 & 25001 \\ \hline
  \end{tabular}  
\end{table}

The size of each subsystem $s \in [S]$ can be measured by $m_s$ and $n_s$ which are the number of rows and columns of $\Ab_s$ in \eqref{OurModel}. 
Since the number of constraints and variables of each subproblem depends on the number of phases, the instance like IEEE 8500 that has smaller number of phases at each component of the network has smaller size of subproblems as reported in Table \ref{tab:size_distributed}.

\begin{table}[!h]
\centering
\caption{The size of subsystems measured by $m_s$ and $n_s$ for $s \in [S]$}
\label{tab:size_distributed}
\begin{tabular}{l|cccc}
\hline
                        & IEEE13 &  IEEE37 & IEEE123 & IEEE8500 \\ \hline 
Min $(m_1,\ldots,m_S)$ & $4$  & $3$& $2$& $2$\\ 
Max $(m_1,\ldots,m_S)$ & $22$ & $42$& $42$& $18$ \\
Mean $(m_1,\ldots,m_S)$ & $9.08$&  $12.43$& $7.34$ & $3.44$ \\ 
Stdev $(m_1,\ldots,m_S)$ & $4.42$&  $8.67$& $4.43$ & $2.66$\\ 
Sum $(m_1,\ldots,m_S)$ & $453$&  $1206$& $1834$ & $86108$\\  \hline
Min $(n_1,\ldots,n_S)$ & $8$  & $6$& $4$& $4$\\ 
Max $(n_1,\ldots,n_S)$ & $34$ & $48$& $57$&$24$ \\
Mean $(n_1,\ldots,n_S)$ & $16.1$&  $20.22$& $13.16$ & $6.69$\\ 
Stdev $(n_1,\ldots,n_S)$ & $5.14$&  $7.35$& $6.5$ & $3.21$\\ 
Sum $(n_1,\ldots,n_S)$ & $805$ & $1961$& $3289$ & $167394$\\     
\hline
\end{tabular}  
\end{table}

As a default setting for Algorithm \ref{alg:admm}, we set $\epsilon_{\text{rel}}=10^{-3}$ in \eqref{termination}, and $\rho=100$.
As an initial point to start (in line 1 of Algorithm \ref{alg:admm}), for each $s \in [S]$, we set $\lambdab^1_s \gets 0$ and set each element of $\xb^1_s$ as either (i) zero if it does not have bounds, (ii) the average value of lower and upper bounds if it has bounds, and (iii) one if it is related to the voltage magnitude.

\subsection{Primal and dual residuals}
\begin{figure}[h!]
\centering    
\includegraphics[width=.45\linewidth]{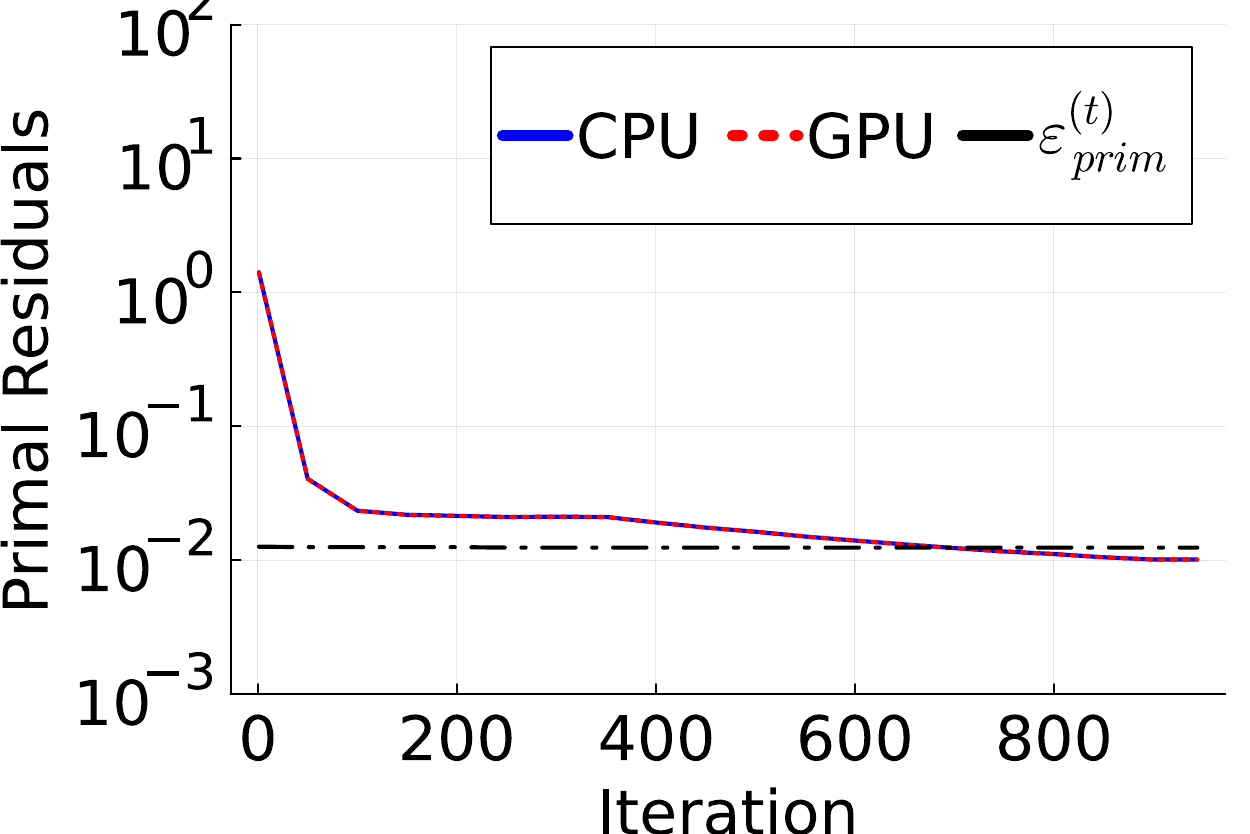}    
\includegraphics[width=.45\linewidth]{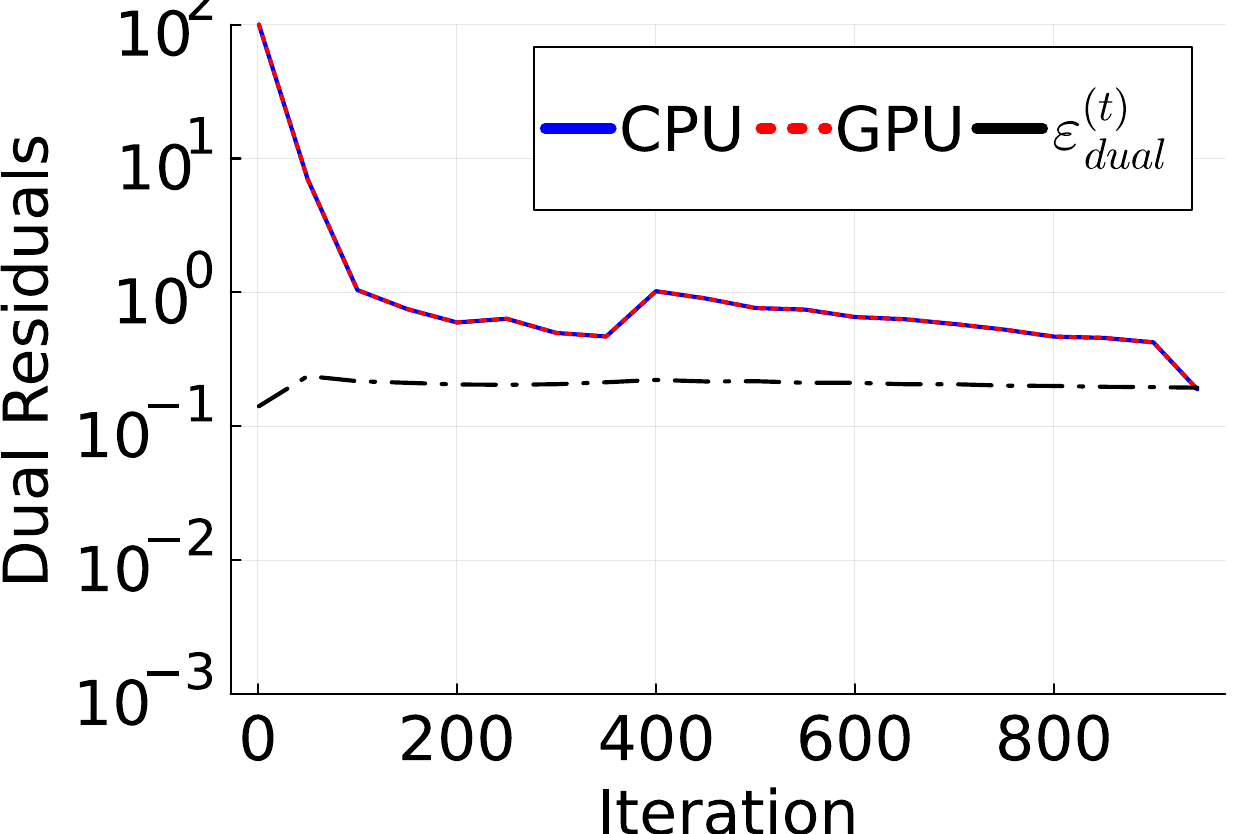}    
\caption{Primal (left) and dual (right) residuals at each iteration of Algorithm \ref{alg:admm} when a CPU and a GPU, respectively, is used for solving the IEEE 13 instance.}
\label{fig:residuals}
\end{figure}

In this section we report the primal and dual residuals as described in \eqref{termination} computed at each iteration of Algorithm \ref{alg:admm} when a CPU and a GPU are used for computation.
Note that a double-precision is used for the CPU computation while a single-precision is used for the GPU computation.
Despite of the difference on the precision, the primal and dual residuals produced by CPU and GPU are close to each other as reported in Figure \ref{fig:residuals}.
Therefore, the total number of iterations required for a convergence to an optimal solution is the same for both CPU and GPU computation, which is reported in Table \ref{tab:comparison}.

\begin{table}[!h]
\centering
\caption{Total number of iterations.}
\label{tab:comparison}
\begin{tabular}{l c c c c}
\hline
& IEEE13& 	IEEE37& 	IEEE123& 	IEEE8500 \\ \hline
Iteration	& 944	& 1774	& 3496	& 15817 \\
\hline
\end{tabular}  
\end{table} 

\begin{figure*}[!h]
\centering
\begin{subfigure}[b]{0.24\textwidth}
\centering
\includegraphics[width=\textwidth]{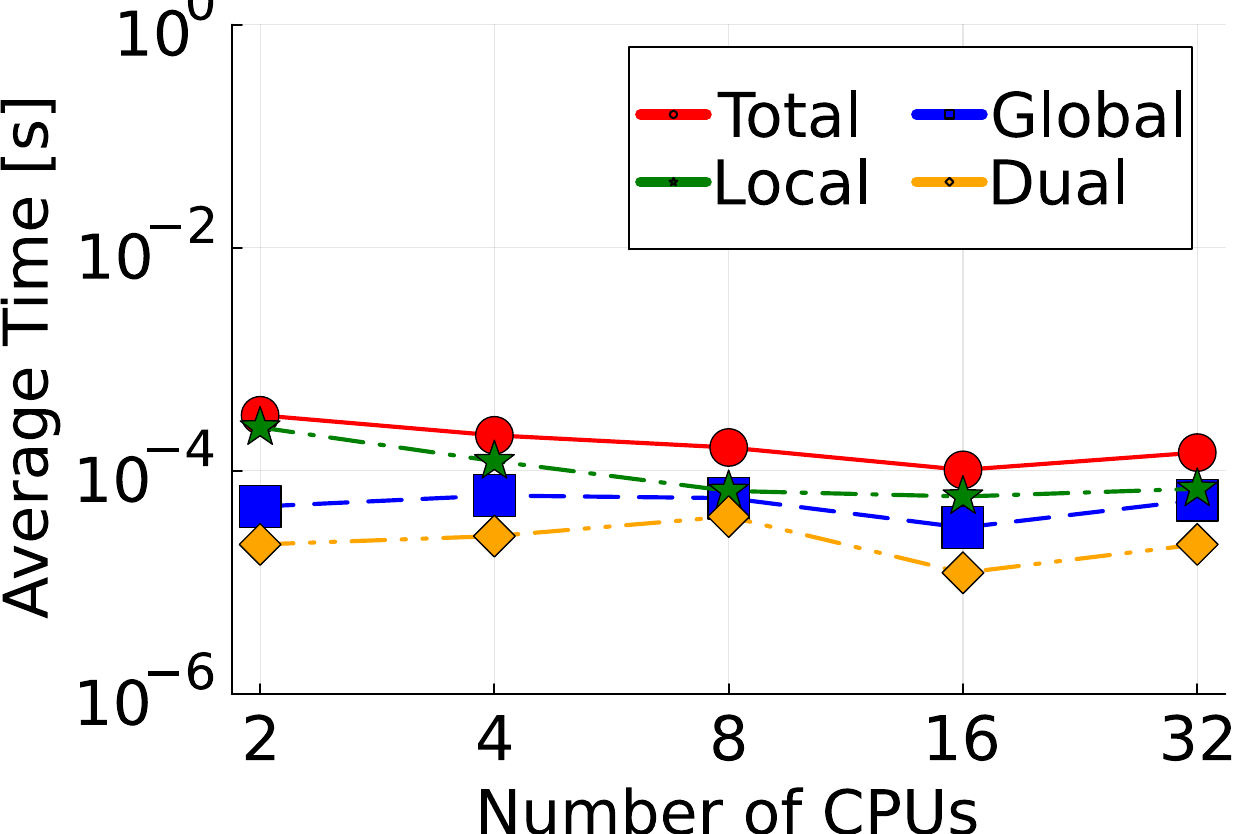}  
\includegraphics[width=\textwidth]{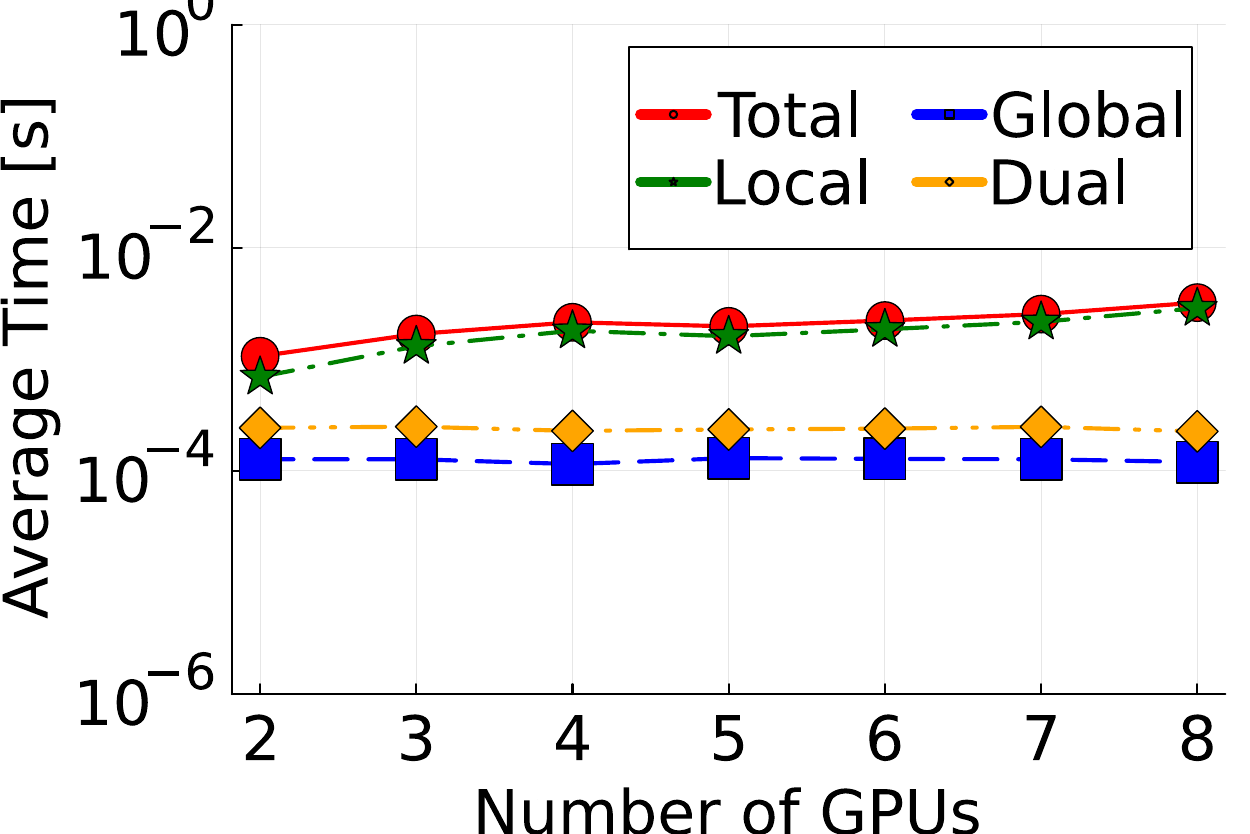}   
\includegraphics[width=\textwidth]{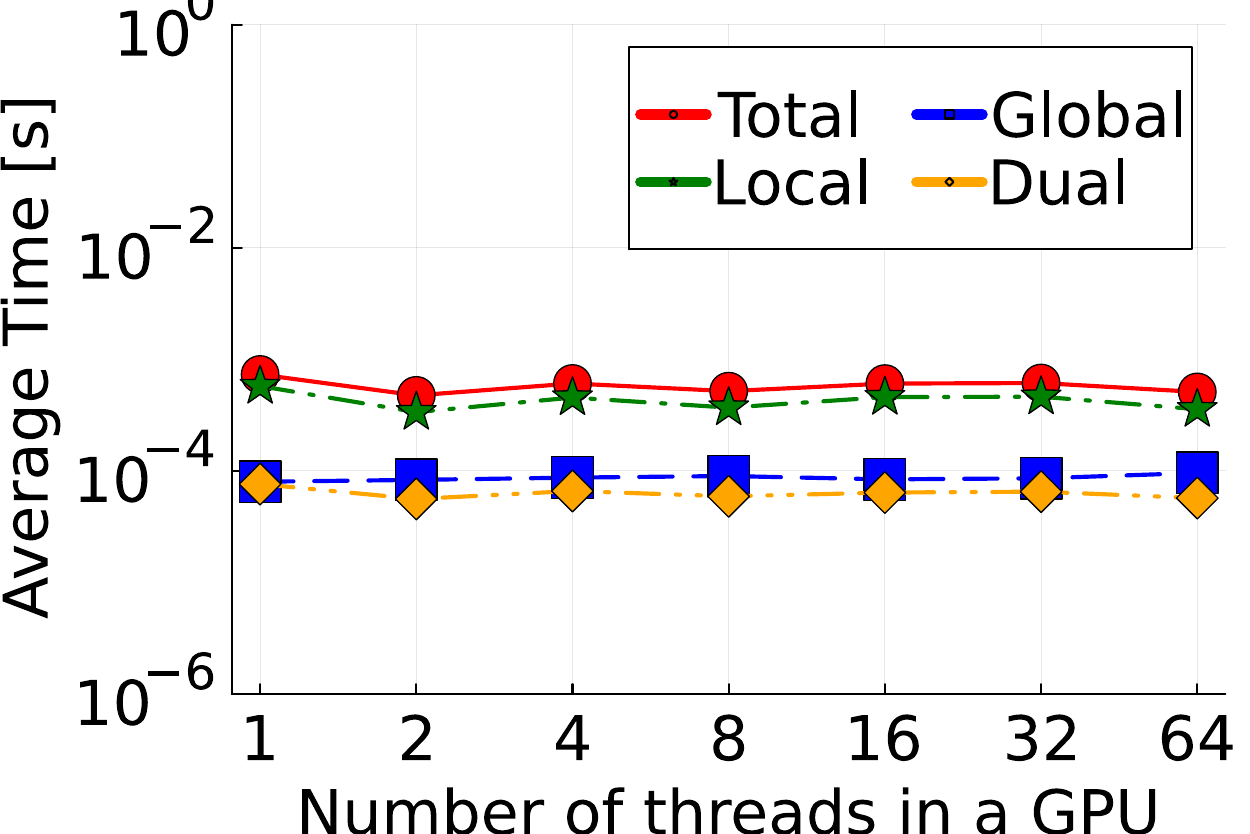}  
\caption{IEEE13}  
\end{subfigure}
\hfill
\begin{subfigure}[b]{0.24\textwidth}
\centering
\includegraphics[width=\textwidth]{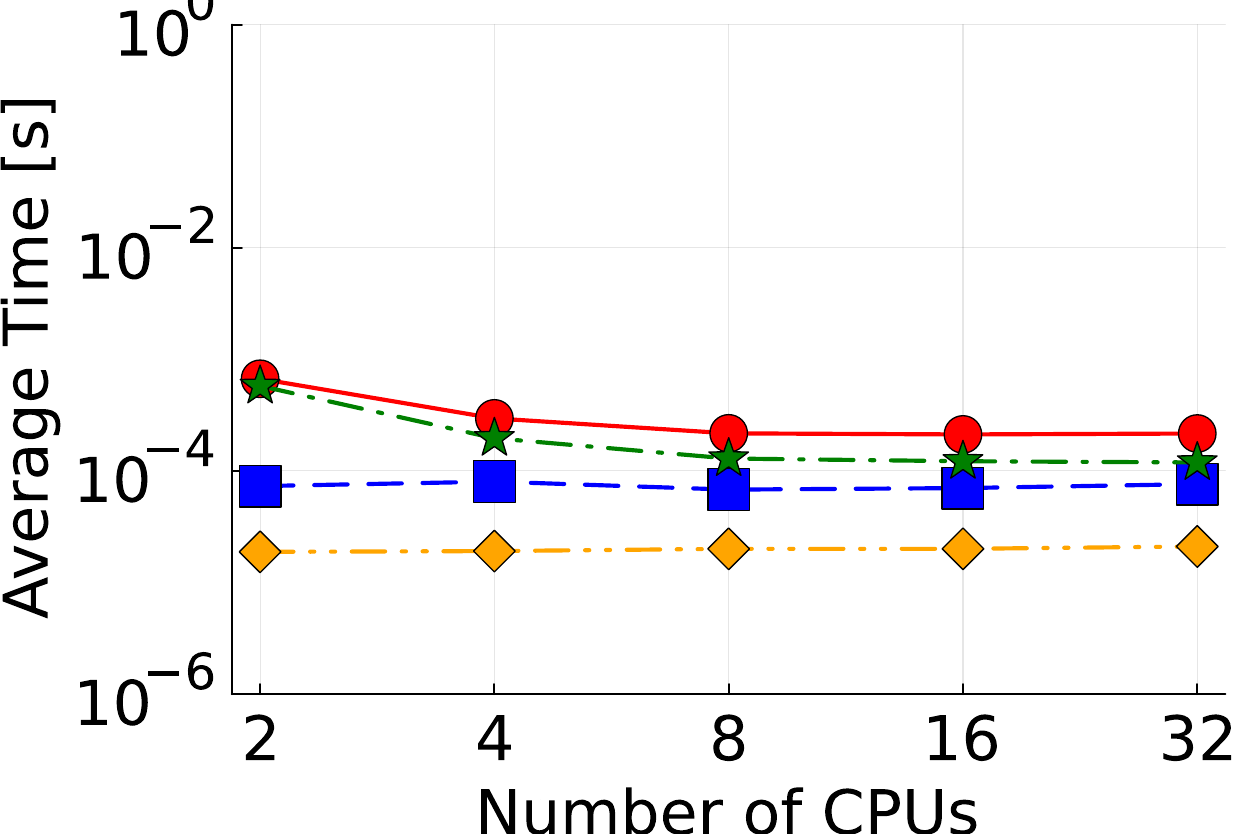} 
\includegraphics[width=\textwidth]{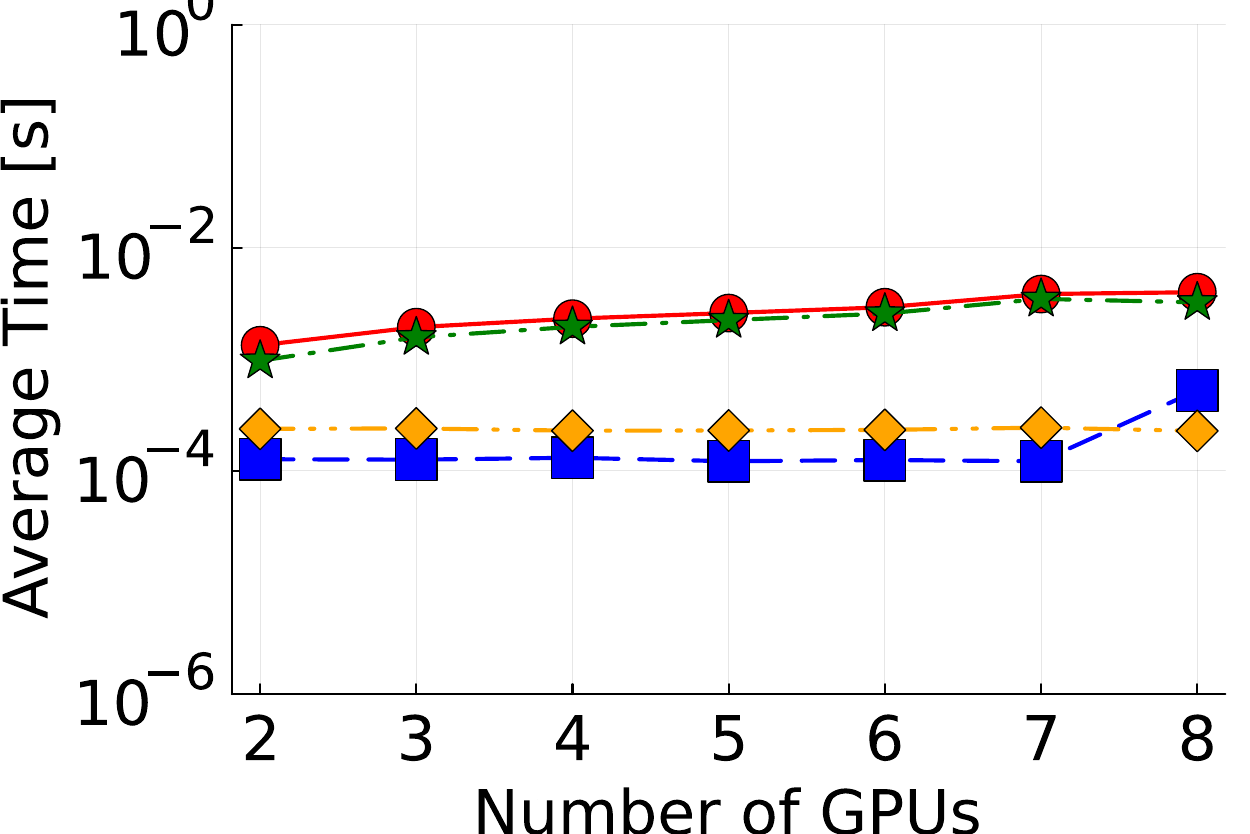}  
\includegraphics[width=\textwidth]{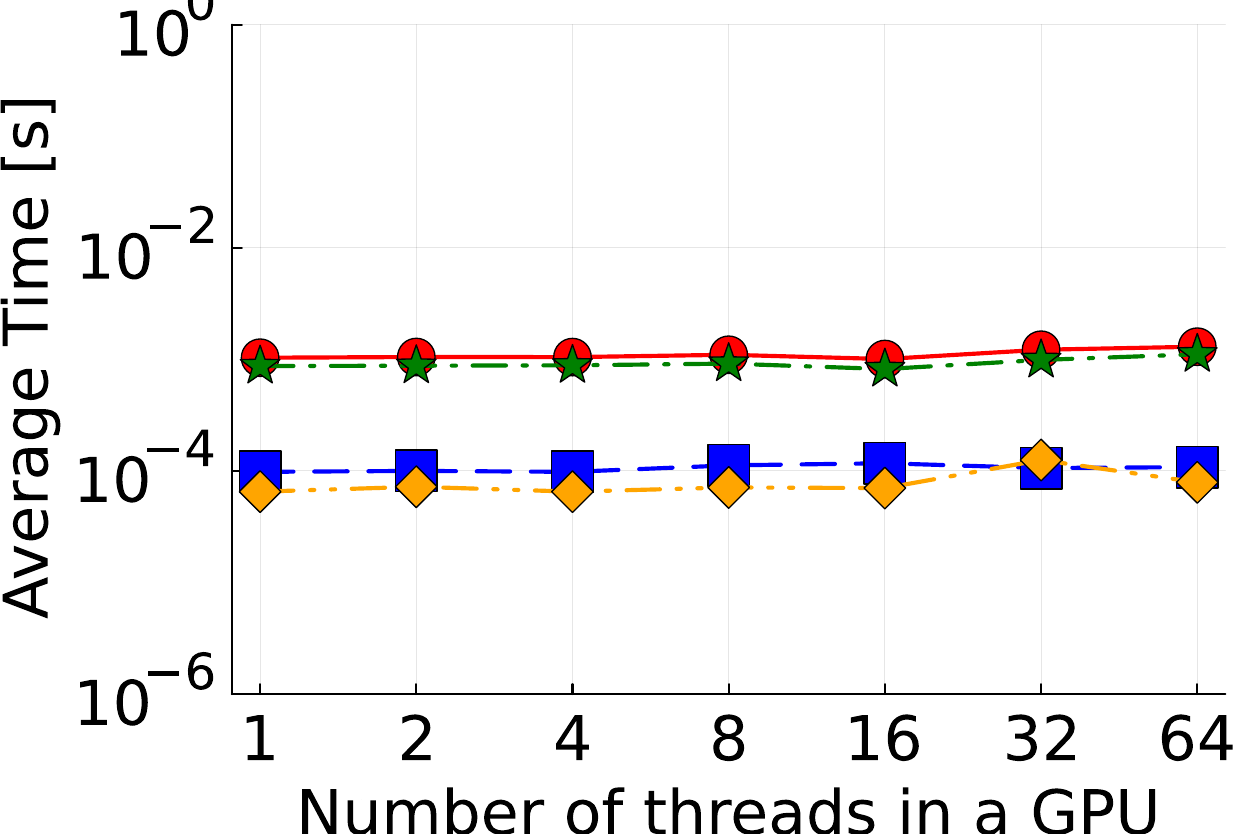}  
\caption{IEEE37}  
\end{subfigure}
\begin{subfigure}[b]{0.24\textwidth}
\centering
\includegraphics[width=\textwidth]{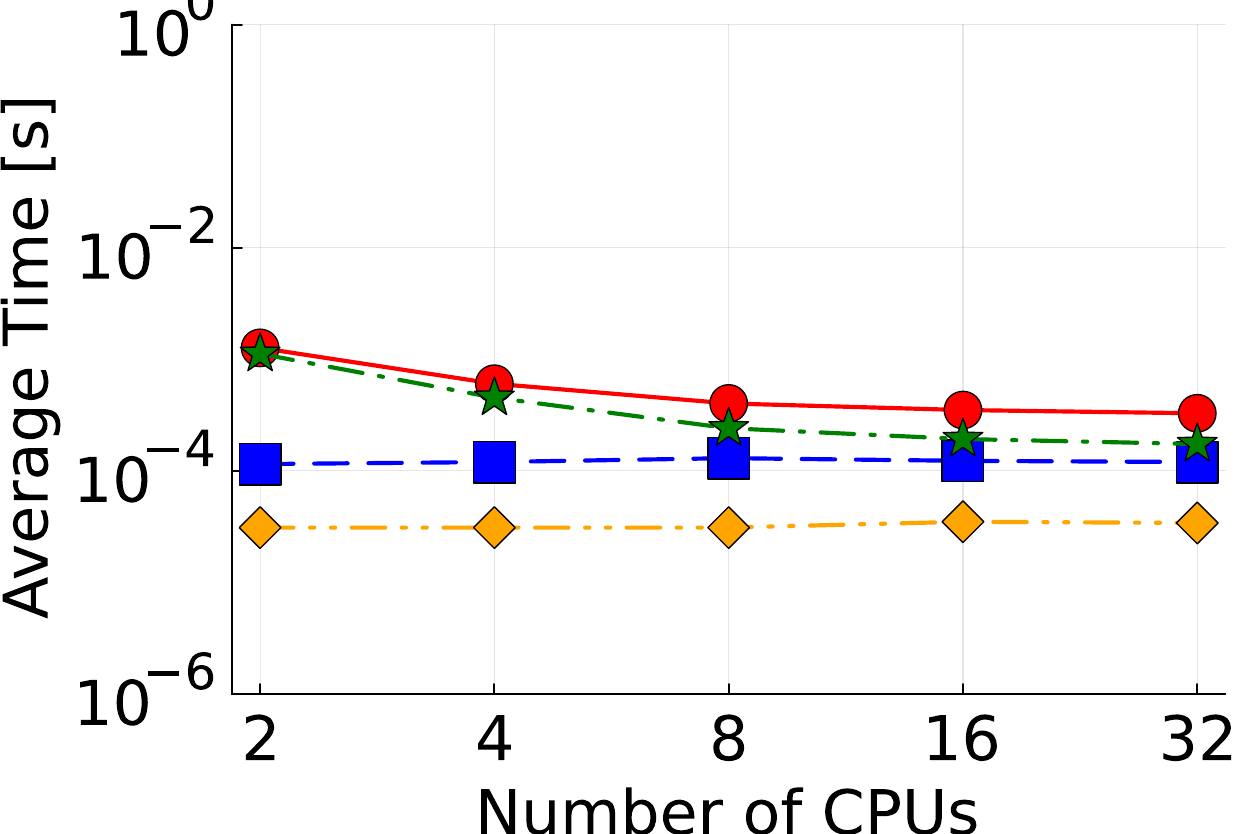} 
\includegraphics[width=\textwidth]{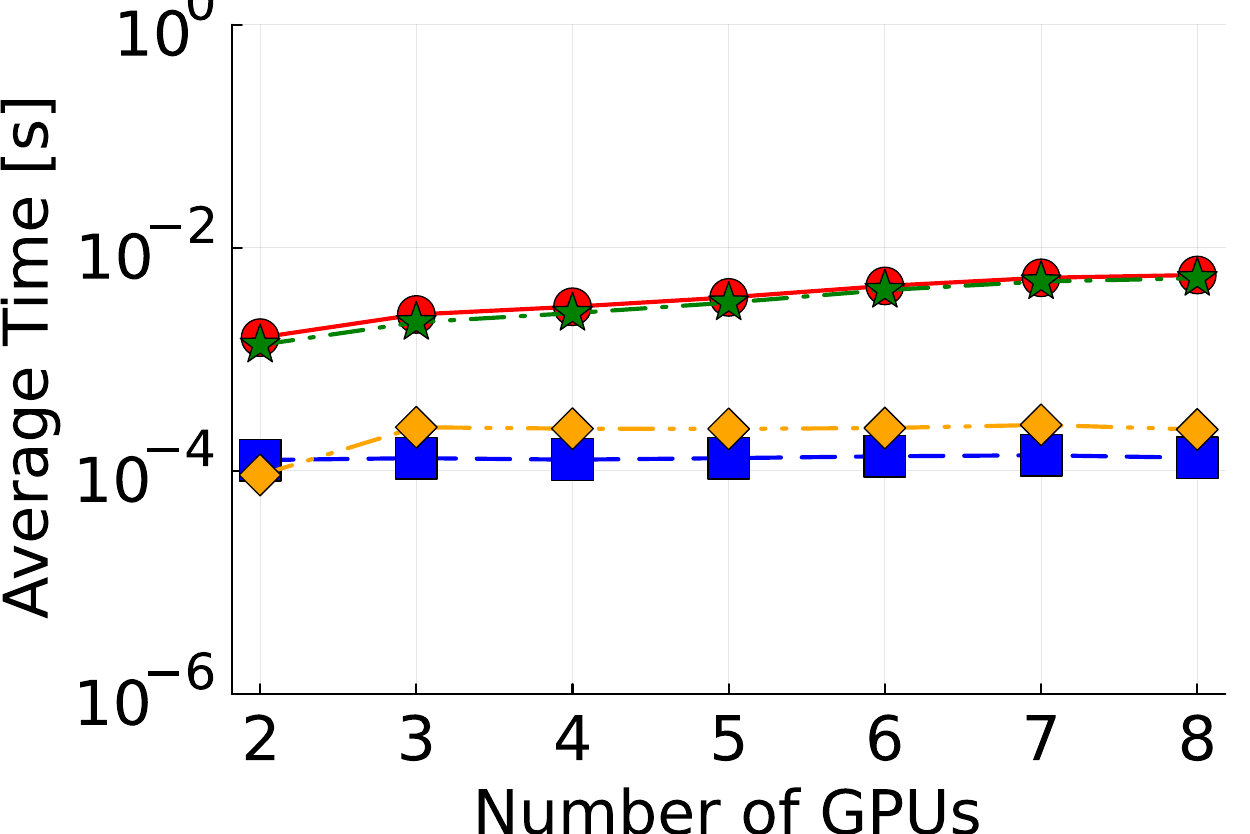}  
\includegraphics[width=\textwidth]{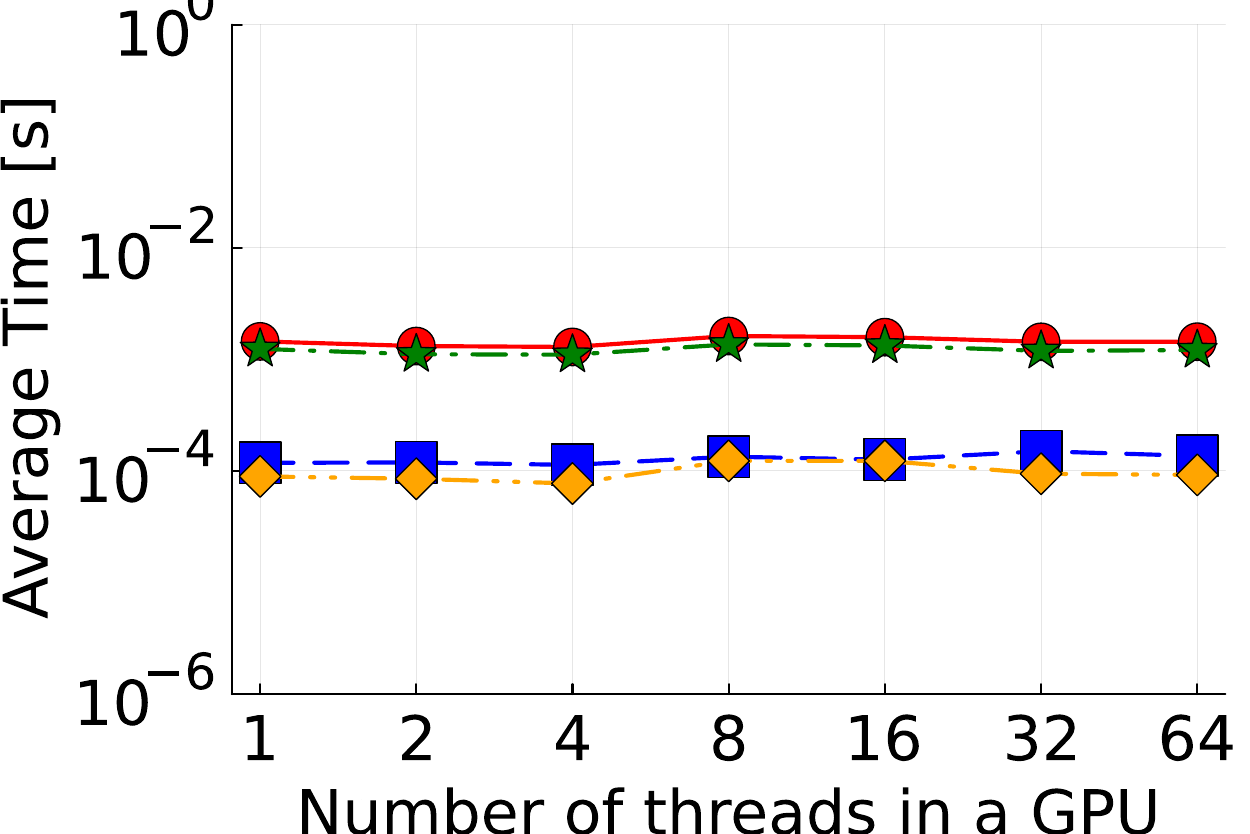}  
\caption{IEEE123}  
\end{subfigure}
\begin{subfigure}[b]{0.24\textwidth}
\centering  
\includegraphics[width=\textwidth]{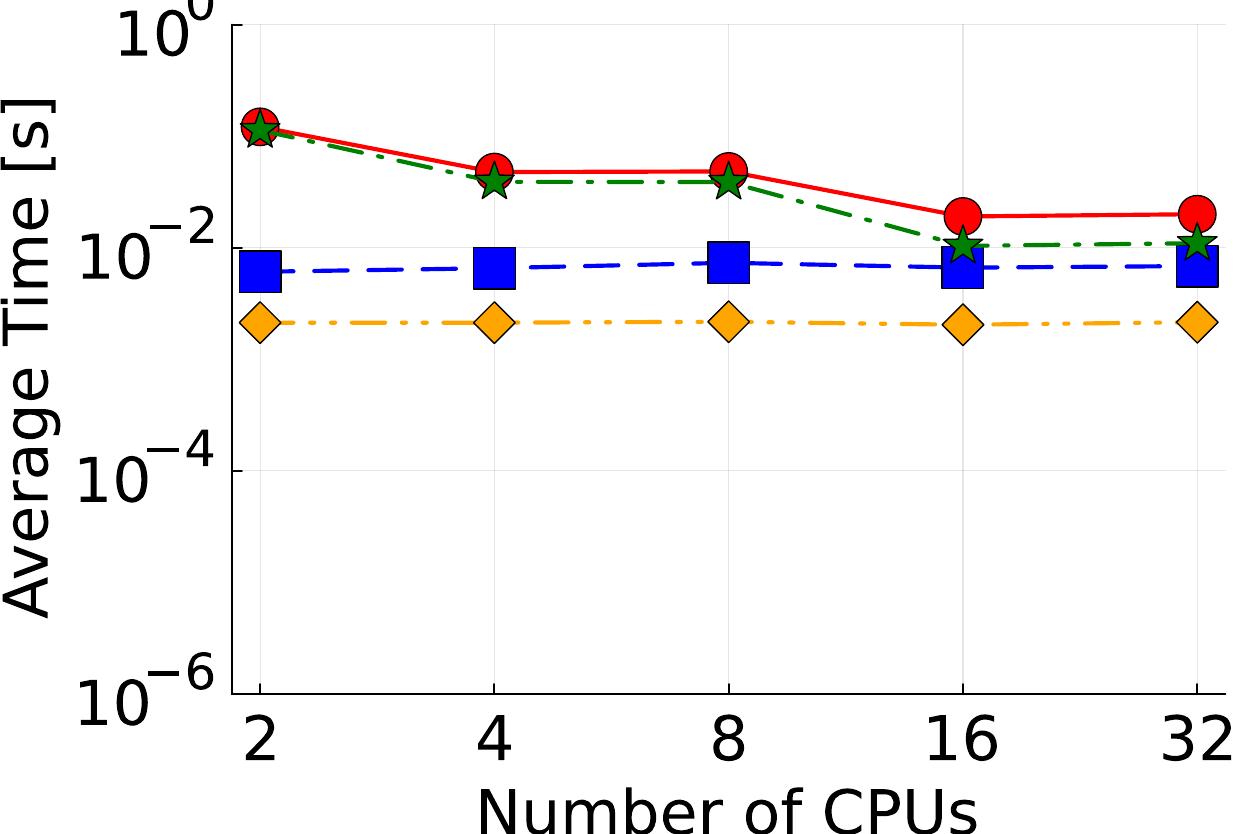} 
\includegraphics[width=\textwidth]{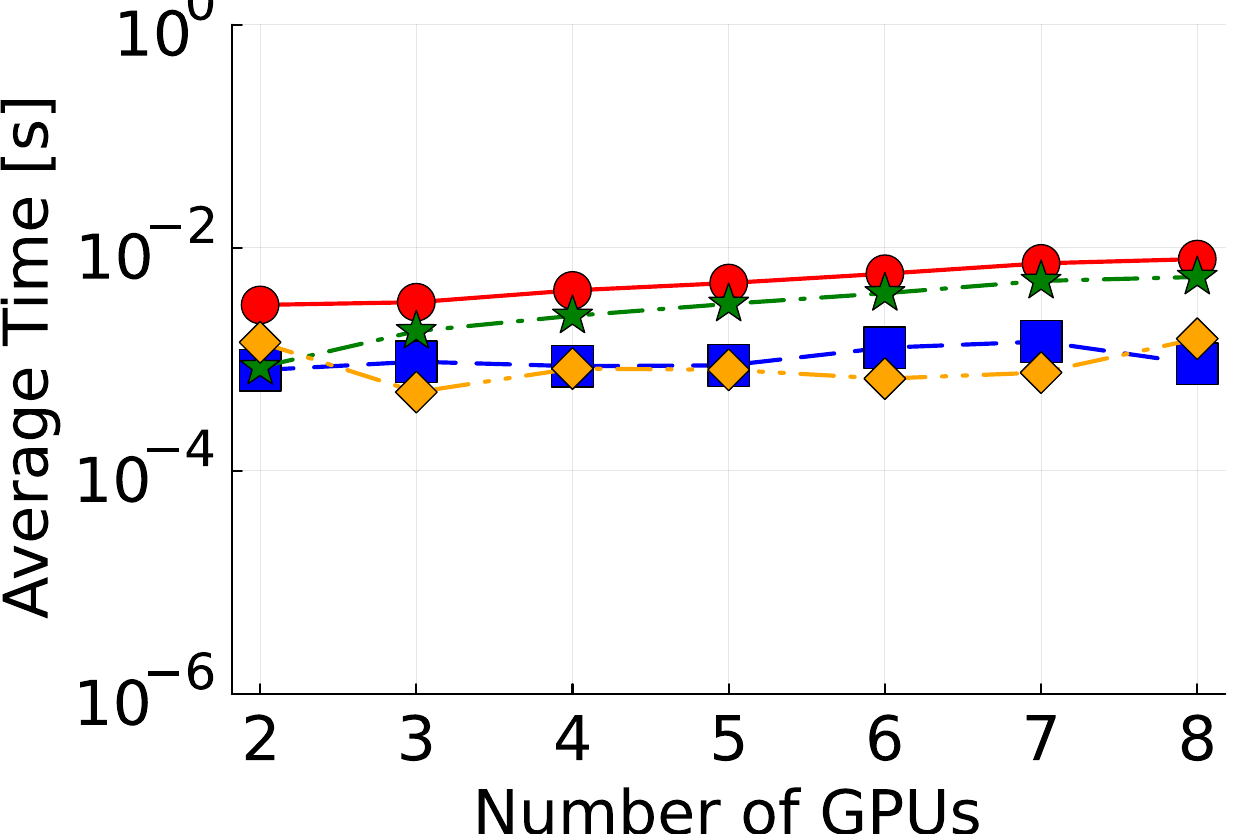} 
\includegraphics[width=\textwidth]{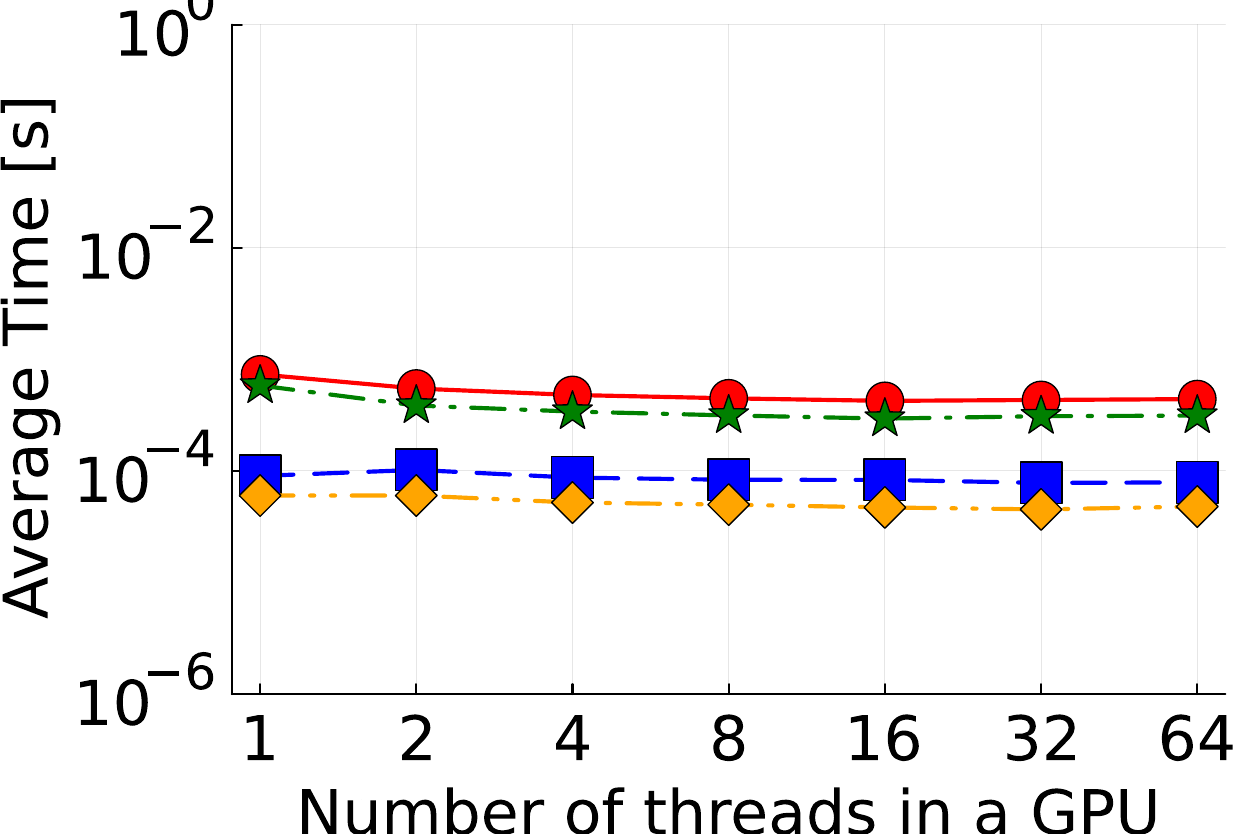} 
\caption{IEEE8500}  
\end{subfigure}          
\caption{Average time for conducting global, local, and dual updates per iteration and their summation (referred as total) when multiple CPUs (top), GPUs (middle), and threads in a GPU (bottom), respectively.}
\label{Figure:cpus_gpus}
\end{figure*}

\subsection{Performance Comparisons}
Given the equivalence in the number of iterations required for Algorithm \ref{alg:admm} to converge to an optimal solution when executed on both CPU and GPU platforms, it becomes imperative to enhance the efficiency of computation per iteration. In pursuit of this objective, we conduct parallel computation employing multiple CPUs and GPUs, respectively.
 
In the first row of Figure \ref{Figure:cpus_gpus}, we present the average computation times per iteration for executing the global update \eqref{closed_1} (indicated by blue squares), the local update \eqref{closed_2} (depicted by green stars), and the dual update \eqref{ADMM-3} (represented by yellow diamonds). Subsequently, we compute the average total time per iteration (denoted by red circles) by aggregating these three values. It is noteworthy that with an increasing number of CPUs, we observe a decrease in the time required for the local update, while the times for global and dual updates remain relatively constant. This behavior aligns with expectations, as the augmentation of CPUs for parallel computation expedites the local update process. Nevertheless, it is evident that the scalability of CPU-based computation is constrained by the network's size.

To provide a comparative perspective, in the subsequent row of Figure \ref{Figure:cpus_gpus}, we present the average computation times when multiple GPUs are employed for parallel computation. As the number of GPUs increases, we observe a slight increase in the times for local updates, while those for global and dual updates remain stable. This phenomenon is primarily attributed to the utilization of MPI for communication within the GPU computing cluster, necessitating data transfer from GPU to CPU. In practice, as elucidated in Section \ref{sec:implementation}, the deployment of MPI may be infeasible due to the physical separation of CPUs and GPUs in distinct clusters. In such cases, remote procedure calls, for instance, tRPC, may be employed for communication, resulting in equivalent communication times for both CPUs and GPUs.

For a more comprehensive evaluation of GPU-based computation, excluding the communication overhead, we exclusively employed a single GPU with multiple threads, as reported in the final row of Figure \ref{Figure:cpus_gpus}. This configuration exhibited notably superior performance compared to parallel CPU computation and exhibited scalability with the network's size.

\begin{figure}[!h]
\centering        
\includegraphics[width=0.7\linewidth]{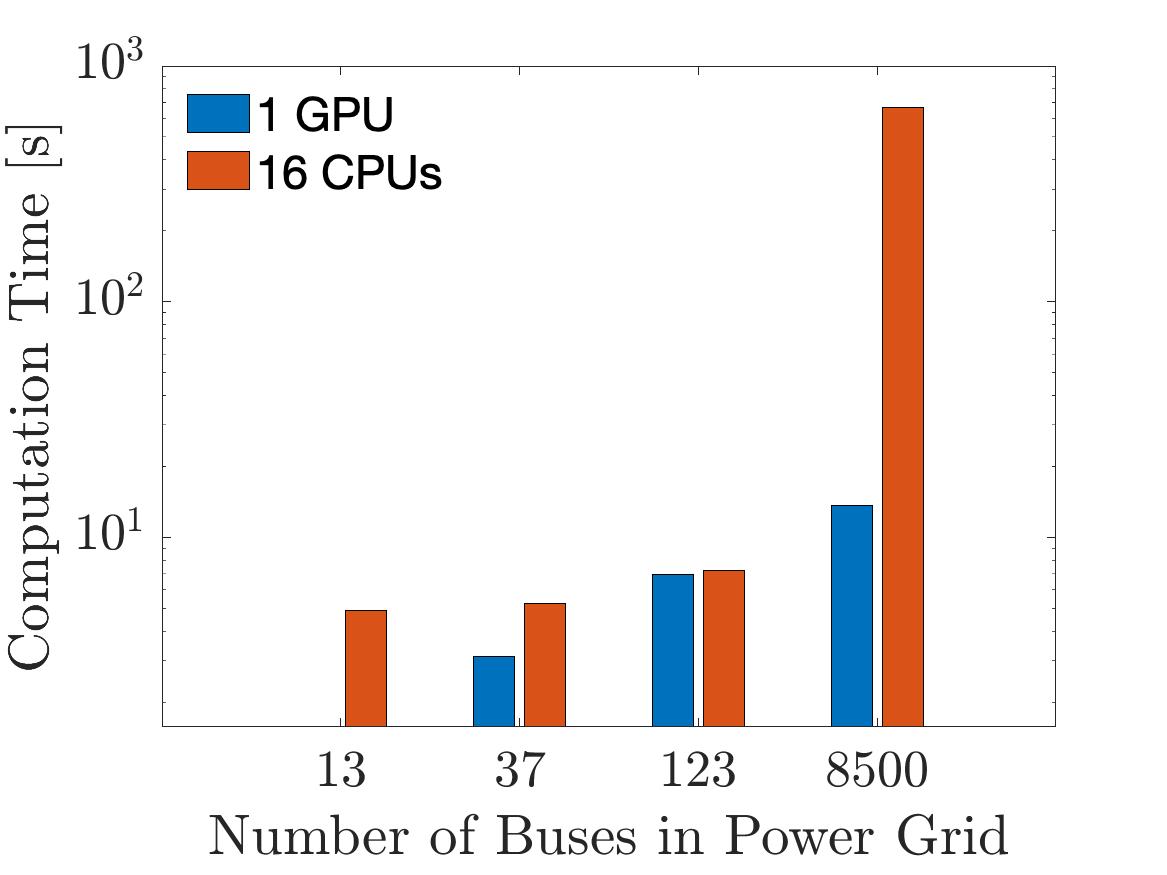}
\caption{Comparison of computation wall-clock time spent by a GPU and 16 CPUs when solving the four instances with the proposed algorithm.
}
\label{fig:time}
\end{figure}

In summary, as presented in Figure \ref{fig:time}, the utilization of a GPU for computation has substantially improved the total computation time, manifesting a fifty-fold enhancement for the IEEE 8500 instance.

\section{Conclusion and Future Work}
We have introduced a GPU-based distributed optimization algorithm tailored for solving the linearized OPF problem within distribution systems. This algorithm builds upon the principles of ADMM and is constructed with subproblems characterized by closed-form expressions relying on matrix operations that can be executed with remarkable efficiency using GPUs.
Notably, the implementation of our approach for GPU has resulted in a substantial reduction in computation time per iteration of the algorithm, achieving a performance gain of up to 50-fold when compared to its CPU-based counterpart, particularly evident when applied to large-scale test instances.
Moreover, the speedup achieved by GPU acceleration would be significantly increasing with much larger instances (e.g., 70,000-bus system as shown in \cite{kim2022accelerated}).

As part of our future research endeavors, we envision the development of a GPU-based distributed optimization algorithm specifically tailored for the second-order conic relaxation of the OPF model. This extension aims to further advance the capabilities of GPU-accelerated optimization techniques in the context of power system analysis and control, potentially opening new avenues for optimizing complex power system problems with enhanced efficiency and scalability.
GPU optimization solvers will also enable the seamless integration of various machine learning techniques (e.g., \cite{zeng2022reinforcement}) that are mainly run on GPUs.

\bibliographystyle{IEEEtran}
\bibliography{reference.bib}

\end{document}